\documentclass[twoside,12pt]{article}
\usepackage{fullpage}
\usepackage{graphicx}
\usepackage{amsmath,amsthm}
\usepackage[mathscr]{eucal}
\usepackage{amssymb,amsfonts,amsbsy}
\usepackage{dsfont,eufrak,mathrsfs}
\usepackage[utf8]{inputenc}
\usepackage{color}
\usepackage[T1]{fontenc}
\usepackage{tikz}
\usepackage{tikz-cd}

\newtheorem{theorem}{Theorem}
\newtheorem{definition}[theorem]{Definition}
\newtheorem{lemma}[theorem]{Lemma}
\newtheorem{proposition}[theorem]{Proposition}
\newtheorem{remark}[theorem]{Remark}
\newtheorem{corollary}[theorem]{Corollary}

\newcommand*{\dsty}{\displaystyle}
\newcommand*{\fun}[3]{#1: #2 \rightarrow #3}

\newcommand*{\ccomp}[1]{\mathcal{C}(#1)}

\title{Isomorphisms between curve graphs of infinite-type surfaces are geometric}
\author{Jes\'{u}s Hern\'{a}ndez Hern\'{a}ndez, Israel Morales, Ferr\'{a}n Valdez}
\date{}

\begin{document}
\maketitle
\begin{abstract}
Let $\phi:\mathcal{C}(S)\to\mathcal{C}(S')$ be a simplicial isomorphism between curve graphs of infinite-type surfaces. In this paper we show that in this situation $S$ and $S'$ are homeomorphic and $\phi$ is induced by a homeomorphism $h:S\to S'$.
\end{abstract}

\section{Introduction}

This is the last of three papers on which the authors study the natural action of the extended mapping class group ${\rm Mod}^*(S)$ of an infinite-type\footnote{A surface is of \emph{infinite type} if its fundamental group is not finitely generated.} surface $S$ on the curve graph $\mathcal{C}(S)$  and whether any isomorphism between curve graphs actually comes from a homeomorphism, see \cite{HHV}, \cite{HHMV}. Our main result is the following:  


\begin{theorem}
	\label{THM:Main}
	Let $S$ and $S'$ be infinite-type connected orientable surfaces with empty boundary and $\phi:\mathcal{C}(S)\to\mathcal{C}(S')$ a simplicial isomorphism. Then $S$ is homeomorphic to $S'$ and $\phi$ is induced by a homeomorphism $h:S\to S'$.
\end{theorem}

As an immediate consequence of this result and Corollary 1.2 in \cite{HHMV} we obtain an analogue for infinite-type surfaces of a foundational well-known result by Ivanov (see Theorem 1 in \cite{Ivanov}):

\begin{theorem}
	\label{THM:surjectivity}
Let $S$ be an infinite-type connected orientable surface with empty boundary. Then every automorphism of the curve graph          $\mathcal{C}(S)$ is induced by a homeomorphism. More precisely, the natural map:
$$
\Psi:{\rm Mod}^*(S)\to {\rm Aut}(\mathcal{C}(S))
$$
is an isomorphism.
\end{theorem}	
	
	It is important to remark that both of these results were known to be true for infinite-type surfaces for which all topological ends carry (infinite) genus \cite{HHV}. With this in mind we highlight the \textbf{main contribution} of this text:	
A new and very simple proof to the fact that every automorphism of the curve graph of an infinite-type surface is geometric. The technology we present is based on \emph{principal exhaustions}, which were introduced in \cite{HHMV}. Roughly speaking, these are nested sequences of finite-type subsurfaces whose union is $S$ and that allow us to ``pull'' classical results about simplicial actions of mapping class groups to the realm of infinite-type surfaces. In particular, the proof that we present of Theorem \ref{THM:Main} makes no use of Dehn-Thurston coordinates and is hence completely independent (and simpler) from the proofs presented in \cite{HHV}.



\emph{Applications}. Given that the curve graph of an infinite-type surface has diameter 2, the natural action of ${\rm Mod}^*(S)$ on $\mathcal{C}(S)$ gives no large scale information. However, the results that we present in this text have found   the following non-trivial applications:
\begin{enumerate}
\item Let $S$ be an infinite-type surface whose genus is finite and at least 4 (\emph{e.g.} a closed surface of genus $5$ to which we have removed a Cantor set) and denote by ${\rm PMod}(S)$ the subgroup of ${\rm Mod}(S)$ consisting of all orientation preserving mapping classes acting trivially on the topological ends of $S$. This group is called \emph{the pure mapping class group} of $S$.
There is a natural monomorphism from ${\rm Mod}^*(S)$ to ${\rm Aut(PMod(S))}$ given by the action of ${\rm Mod}^*(S)$ on the pure mapping class group by conjugation. The following result, obtained by Patel and Vlamis in \cite{PatelVlamis}, is a generalization of a famous result of Ivanov for mapping class groups of surfaces of finite type (see \cite{Ivanov88}): 
\begin{theorem}
	\label{THM:VlamisPatel}
If $S$ is an infinite-type surface of genus at least 4, then the natural monomorphism from ${\rm Mod}^*(S)$ modulo its center to ${\rm Aut(PMod(S))}$ is an isomorphism. 
\end{theorem}
Indeed, Patel and Vlamis show that any automorphism of 
${\rm PMod}(S)$ preserves Dehn twists. As they remark, if every automorphism of $\mathcal{C}(S)$ is induced by a homeomorphism of the surface, then Theorem \ref{THM:VlamisPatel} follows by a standard argument that can be found in Ivanov's original paper [\emph{Ibid}.].

\item Let $\Sigma_g$ be the surface that results from removing a Cantor set from a closed orientable surface of genus $g\geq 0$. In \cite{AramayonaFunar} Aramayona and Funar study $\mathcal{B}_g$ the \textit{asymptotically rigid mapping class group of} $\Sigma_g$. This is a finitely presented subgroup of ${\rm Mod}^*(\Sigma_g)$ which contains the mapping class group of every surface of genus $g$ with nonempty boundary. Using Theorem \ref{THM:surjectivity}, they prove that $\mathcal{B}_g$ is \emph{rigid}, that is: 
\begin{theorem}
	\label{THM:AramayonaFunar}
For every $g<\infty$, ${\rm Aut}(\mathcal{B}_g)$ coincides with the normalizer of $\mathcal{B}_g$ within ${\rm Mod}^*(\Sigma_g)$. 
\end{theorem}
More precisely, their result uses the following lemma, which is proven using Theorem \ref{THM:surjectivity}: 
\begin{lemma}
	\label{Lemma:AramayonaFunar}
If $g<\infty$ then ${\rm Aut(PMod_c}(\Sigma_g))={\rm Mod}^*(\Sigma_g)$.
\end{lemma}
Here ${\rm PMod_c}(\Sigma_g)$ denotes the subgroup of ${\rm Mod}^*(\Sigma_g)$ formed by \emph{compactly supported pure} mapping classes. For details we refer the reader to [\emph{Ibid}.].
\end{enumerate}

It is in the light of these facts that we conjecture that further applications of Theorems \ref{THM:Main} and \ref{THM:surjectivity} to the study of big mapping class groups should exist. In particular, it is a natural to wonder which of the classical applications of Ivanov's theorem have an analog in the realm of infinite-type surfaces.

\textbf{Reader's guide.} In \S \ref{section:preliminaries} we make a short discussion on the general aspects of 
pants decompositions. We also recall the notion of principal exhaustion. In \S \ref{SEC:TopologicalRigidity} we prove several topological properties that are preserved under isomorphisms of curve graphs. 
Finally, in \S \ref{Section:Surjectivity} we prove that these isomorphisms are geometric.

\textbf{Acknowledgements.} The first author was supported during the creation of this article by the UNAM Post-Doctoral Scholarship Program 2016 at the Centro de Ciencias Matem\'{a}ticas, UNAM. The second
author was generously supported by CONACYT and CCM-UNAM. The third author would like to dedicate this work to Marion and Anika Yunu\'en.

\section{Preliminaries}
	\label{section:preliminaries}
The main tools that we use in this text to prove our theorem are pants decompositions and a special kind of exhaustion for infinite-type surfaces called \emph{principal exhaustions}. In what follows we recall the definitions of these objects and the main properties that we need.

Abusing language and notation, we call \textit{curve}, a topological embedding $\mathbb{S}^1\hookrightarrow S$ , the
isotopy class of this embedding and its image on $S$. A curve is said to be \textit{essential} if it is neither homotopic to a point nor to the boundary of a neighbourhood of a puncture. Hereafter all curves considered are essential unless otherwise stated.

A collection of essential curves $\mathcal{L}$ in $S$ is {\it{locally finite}} if for every $x\in S$ there exists a neighbourhood $U_x$ of $x$ which intersects finitely many elements of $\mathcal{L}$. A locally finite collection of pairwise disjoint non-isotopic essential curves is called a \textit{multicurve}. 

\begin{definition} [\sc The Curve Graph]
The curve graph of $S$, denoted by $\mathcal{C}(S)$, is the \emph{abstract} graph whose vertices are 
isotopy classes of essential curves in $S$, and 
two vertices $\alpha$ and $\beta$ span an edge if the corresponding curves are disjoint modulo homotopy.  We denote the set of vertices of $\mathcal{C}(S)$ by $\mathcal{C}^0(S)$.
\end{definition}

\begin{remark}
The curve graph is the 1--skeleton of the \emph{curve complex}, that is, the abstract simplicial complex whose simplices are multicurves of finite cardinality. The curve complex is a \emph{flag} complex, in particular it is completely determined by its 1-skeleton and for this reason in this text we restrain our discussion to the curve graph.
\end{remark}


\begin{definition}[\sc Pants decomposition]

	\label{DEF:PPD}
A maximal (w.r.t. inclusion) multicurve  is called a \emph{pants decomposition}.		
\end{definition}
In this text we call both a maximal multicurve $P=\{\alpha_k\}_{k\in K}$ and its image $\{[\alpha_k]\}_{k\in K}$ in $\mathcal{C}^0(S)$, a  pants decomposition. The following lemma gives a simplicial characterization of pants decompositions.

\begin{lemma}\label{CarPP}
Let $S$ be a surface and $P=\{a_k\}_{k\in\mathbb{N}}\subset \mathcal{C}^0(S)$. Then $P=\{a_k\}_{k\in\mathbb{N}}$ is pants decomposition for $S$ if and only if this collection satisfies: 
 \begin{enumerate}
 \item $i(a_k,a_l)=0$ for all $k,l \in \mathbb{N}$,
 \item  for each $a\in  \mathcal{C}^0(S)\backslash P$ we have:
 \begin{enumerate}
 \item $i(a,a_k)\neq 0$ for some $k\in \mathbb{N}$, and
 \item $\vert \{k\in\mathbb{N}:\, i(a,a_k)\neq 0\} \vert<\infty$.
 \end{enumerate}
 
 \end{enumerate}
\end{lemma}

\begin{proof}
The necessity condition is easily verified. We prove the sufficiency. Since $S$ is of infinite type, it is uniformized by the Poincar\'e disc and hence we can choose a complete Riemannian metric on $S$ with negative constant curvature. Since the metric is complete, we can choose for each $a_k$ the only geodesic representative $\alpha_k$ in its class. Conditions 1 and 2.a above assure that $\{\alpha_k\}_{k\in\mathbb{N}}$ is a maximal collection of non-isotopic curves. To see that this collection is locally finite we proceed as follows. Let $s\in S$ be a point and $N\subset S$ a compact finite-type subsurface containing $s$ such that each connected component of $\partial N$ is a closed geodesic that is essential in $S$. Suppose first that there exists finitely many curves in $\{\alpha_k\}_{k\in\mathbb{N}}$ which intersect $\partial N$. In this case one can easily find a neighbourhood $U_s$ of $s$ in $N$ which intersects only finitely many elements of $\{\alpha_k\}_{k\in\mathbb{N}}$. On the other hand, if infinitely many elements of $\{\alpha_k\}_{k\in\mathbb{N}}$ intersect $\partial N$ then we get a contradiction with condition 2.b. 
\end{proof}

\begin{corollary}
	\label{COR:PreDPP}
Let $\phi:\mathcal{C}(S)\to\mathcal{C}(S')$ be a simplicial isomorphism. Then $\phi$ sends pants decompositions in $S$ to pants decompositions in $S'$.
\end{corollary}

Let $\mathcal{L}$ be a multicurve. We say that $\mathcal{L}$ \textit{bounds a subsurface} $\Sigma$ of $S$, if the elements of $\mathcal{L}$ are \textbf{exactly all} the boundary curves of the closure of $\Sigma$ on $S$. Also, we say that $\Sigma$ \textit{is induced by} $\mathcal{L}$ if there exists a subset $\mathcal{M} \subset \mathcal{L}$ that bounds $\Sigma$ \textbf{and} there are no elements of $\mathcal{L} \backslash \mathcal{M}$ in its interior.



\begin{remark}
	\label{RMK:PPD}

Given that every pants decomposition $P$ is locally finite, the index set $K$ in $P=\{\alpha_k\}_{k\in K}$ is at most countable. Moreover, every connected component of $S \backslash P$ is homeomorphic to the thrice-punctured sphere. In particular, every closed subsurface $\Sigma$ induced by $P$ is homeomorphic to the compact surface of genus zero and three boundary components. This topological surface is called a \emph{pair of pants}.
\end{remark}

Recall that an essential curve $\alpha$ is called \emph{separating} if $S \backslash \alpha$ is disconnected. A separating curve $\alpha$ is called \emph{outer} if it bounds a twice-punctured disc.

Let $P$ be a pants decomposition, and let $\alpha, \beta \in P$. We say $\alpha$ and $\beta$ are adjacent w.r.t. $P$ if there exists a subsurface $\Sigma$ induced by $P$ such that $\alpha$ and $\beta$ are two of it boundary curves.

In all the proofs of our main theorems we use the following graph associated to a  pants decomposition.

\begin{definition}[\sc Adjacency graph]
	\label{DEF:AdjacencyGraph}
	Let $S$ be a surface and $P$ a pants decomposition of $S$. We define the {\it{adjacency graph}} of $P$, denoted by $\mathcal{A}(P)$, as the abstract simplicial graph whose set of vertices is $P$ and where two curves $\alpha$ and $\beta$ span an edge if they are adjacent w.r.t. $P$.
\end{definition}

Given a pants decomposition $P$ and a subset $Y\subseteq P$, we denote by $V(Y)$ the set of vertices in $\mathcal{A}(P)$ defined by elements in $Y$.

Finally, we recall a particular way to exhaust infinite-type surfaces that is used to prove that every automorphism of the curve graph is geometric.

\begin{definition}[\sc Principal exhaustion]\label{DEF:PrincipalExhaustion}
Let $\{S_{i}\}_{i\in\mathbb{N}}$ be an (set-theoretical) increasing sequence of \emph{open connected} subsurfaces of $S$. We say $\{S_{i}\}_{i\in\mathbb{N}}$ is a \textit{principal exhaustion of} $S$ if $S = \dsty \bigcup_{i \geq 1} S_{i} $ and for all $i \geq 1$ it satisfies the following conditions:
\begin{enumerate}
 \item $S_{i}$ is a surface of finite topological type,
 \item $S_{i}$ is contained in the interior of $S_{i+1}$,
 \item $\partial S_{i}$ is the finite union of pairwise disjoint essential separating curves on $S$
 \item each connected component of $S_{i+1} \backslash \overline{S_{i}}$ has complexity at least 4, and
 \item each connected component of $S\setminus \overline{S_i}$ is of infinite topological type.
\end{enumerate}
\end{definition}


\section{Topological properties}
	\label{SEC:TopologicalRigidity}

In this section we prove several topological properties preserved under an isomorphism $\phi: \mathcal{C}(S) \to \mathcal{C}(S')$. All surfaces in this section are of infinite type, unless otherwise stated.




The following two propositions and lemma can be deduced from the work of Shackleton \cite{Shack}. More precisely, Propositions \ref{PROP:GraphIsomorph} and \ref{PROP:SameGenus} below are Lemmas 8 and 12 in [\emph{Ibid}.]; on the other hand, Lemma \ref{Lemma:SpecialToSpecial} below follows from Lemmas 9 and 10 in [\emph{Ibid}.] and the fact that $\phi$ is an isomorphism.

As a matter of fact, Shackelton does not work on the context of infinite-type surfaces but the arguments that he uses to prove these results are of local nature and hence can be immediately extrapolated to all infinite-type surfaces. For the sake of completeness, we include a sketch of proof in each case. 
\begin{proposition}
	\label{PROP:GraphIsomorph}
Let $\phi:\mathcal{C}(S)\to\mathcal{C}(S')$ be a simplicial isomorphism between curve graphs of infinite-type surfaces. Then $\phi$ induces a graph isomorphism $$
\tilde{\phi}:\mathcal{A}(P)\to\mathcal{A}(\phi(P))
$$ for any pants decomposition $P$ of S.
\end{proposition}

\emph{Sketch of proof}. Since pants decompositions are maximal multicurves, $\tilde{\phi}$ is a biyective correpondence between the set of vertices of $\mathcal{A}(P)$ and the vertices of $\mathcal{A}(\phi(P))$. Then we only need to check that $\tilde{\phi}$ and $\tilde{\phi}^{-1}$ preserve edges, but this follows from the fact that any two vertices $\alpha$ and $\beta$ are adyacent in $\mathcal{A}(P)$ if and only if there exist a curve $\gamma$ in $S$ that intersects $\alpha$ and $\beta$ but does not intersect any other element in $P\backslash \{\alpha,\beta\}$. 
\qed

\begin{lemma}
	\label{Lemma:SpecialToSpecial}
Let $\phi:\mathcal{C}(S)\to\mathcal{C}(S')$ be a simplicial isomorphism between curve graphs of infinite-type surfaces $S$ and $S'$. Then $\phi$ maps nonouter separating curves to nonouter separating curves, nonseparating curves to nonseparating curves and hence outer curves to outer curves.
\end{lemma}

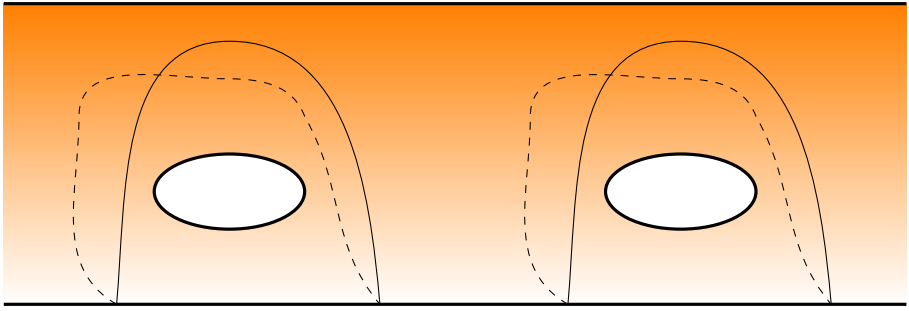
\begin{figure}[!ht]
\begin{center}
\begin{tikzpicture}[scale=2]

\shadedraw[top color=orange, draw=none]
(0,0)--(6,0)--(6,2)--(0,2)--cycle;
\draw[very thick] (0,0)--(6,0);
\draw[very thick] (0,2)--(6,2);
\draw[fill=white, very thick] (1.5,.75) ellipse (.5 and .25);
\draw[fill=white, very thick] (4.5,.75) ellipse (.5 and .25);


\draw (.75,0) to[out=85,in=180] (1.5,1.75) to[out=0,in=95] (2.5,0);
\draw[dashed] (.75,0) to[out=150,in=270] (.5,1.25) to[out=90,in=180] (1.5,1.5) to[out=0,in=110] (2,1.25) to[out=300,in=130] (2.5,0);

\draw[xshift=3cm] (.75,0) to[out=85,in=180] (1.5,1.75) to[out=0,in=95] (2.5,0);
\draw[dashed,xshift=3cm] (.75,0) to[out=150,in=270] (.5,1.25) to[out=90,in=180] (1.5,1.5) to[out=0,in=110] (2,1.25) to[out=300,in=130] (2.5,0);
\end{tikzpicture}
\end{center}
\caption{Curves capturing the genus of the surface.}
\label{FIG:CurvesCatchGenus}
\end{figure}

\emph{Sketch of proof}. Nonouter separating curves are cut vertices of the graph $\mathcal{A}(P)$, for any pants decomposition $P$ and \emph{viceversa}. Since $\tilde{\phi}$ is an isomorphism, cut vertices must go to cut vertices. On the other hand, outer curves are vertices of degree at most 2 in $\mathcal{A}(P)$, for any pants decomposition $P$. Given that $\phi$ is an isomorphism, nonseparating curves cannot be mapped to nonouter separating curves. If for a nonseparating curve $\alpha$, $\phi(\alpha)$ were an outer curve, then  we could find a pants decomposition $P$ containing $\alpha$ for which the vertex corresponding to $\alpha$ has degree four. This contradicts the fact that $\phi$ is an isomorphism. \qed
\begin{proposition}
	\label{PROP:SameGenus}
Let $\phi:\mathcal{C}(S)\to\mathcal{C}(S')$ be a simplicial isomorphism between curve graphs of infinite-type surfaces. Then $S$ and $S'$ have the same genus.
\end{proposition}
\emph{Sketch of proof}.
Let $\mathcal{L}$ be a multicurve in $S$ such that each curve in $\mathcal{L}$ bounds a once-punctured torus in $S$ and $S \backslash \mathcal{L}$ has only one connected component of infinite type and genus zero. In other words, $\mathcal{L}$ is the multicurve that ``captures'' all genus in $S$, see Figure \ref{FIG:CurvesCatchGenus}. Hence the genus of $S$ is equal to the cardinality of $\mathcal{L}$. By Lemma 12 in [\emph{Ibid.}], for each $\alpha\in \mathcal{L}$ the curve $\phi(\alpha)$ bounds a once-punctured torus in $S'$ induced by $\phi(\mathcal{L})$. Hence $genus(S)\leq genus(S')$. As $\phi$ is an isomorphims we obtain the equality.
\qed


Recall that two curves $\{\alpha,\beta\}$ form a \emph{peripheral pair} if they bound a once-punctured annulus.

\begin{figure}[!ht]
\begin{center}
\begin{tikzpicture}[scale=2]
\shadedraw[top color=orange, draw=none] (0,0) to[out=90,in=270] (0,.5) to[out=180,in=270] 
(-.5,1.5) to[out=90,in=270] 
(.5,3) to[out=270,in=90] 
(1.5,1.5) to[out=270,in=0] 
(1,.5) to[out=270,in=90] 
(1,0);
\draw[very thick] (0,0) to[out=90,in=270] (0,.5) to[out=180,in=270] 
(-.5,1.5) to[out=90,in=270] 
(.47,3);

\draw[very thick] (.53,3) to[out=270,in=90] 
(1.5,1.5) to[out=270,in=0] 
(1,.5) to[out=270,in=90] 
(1,0);
\draw[fill=white, very thick] (.5,1.5) ellipse (.25 and .5);
\draw (-.5,1.5) to[out=340,in=200] (.25,1.5);
\draw[dashed] (-.5,1.5) to[out=20,in=160] (.25,1.5);
\draw[xshift=1.25cm] (-.5,1.5) to[out=340,in=200] (.25,1.5);
\draw[xshift=1.25cm, dashed] (-.5,1.5) to[out=20,in=160] (.25,1.5);
\draw (0,.5) to[out=340,in=200] (1,.5);
\draw[dashed] (0,.5) to[out=20,in=160] (1,.5);

\node at (-0.75,1.5) {$\alpha$};
\node at (1.75,1.5) {$\beta$};
\node at (-0.25,.25) {$\gamma$};
\draw[fill=white, thick] (.5,3) ellipse (.05 and .025);
\end{tikzpicture}
\end{center}
\caption{Nonseparating curves $\alpha$ and $\beta$ forming a peripheral pair.}
\label{FIG:NonSepPeriphPair}
\end{figure}
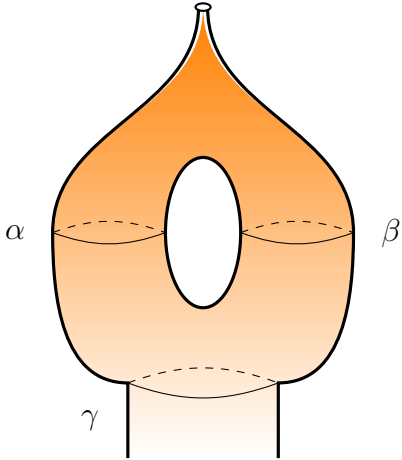

\begin{proposition}
	\label{PROP:PeripheralPairs}
Let $\phi:\mathcal{C}(S)\to\mathcal{C}(S')$ be a simplicial isomorphism between curve graphs of infinite-type surfaces. Then $\phi$ maps peripheral pairs to peripheral pairs. 
\end{proposition}
\begin{proof}
Remark that if $\{\alpha,\beta\}$ is a peripheral pair then both curves forming it have to be either separating or nonseparating. Therefore we only consider the following three cases: (1) both $\alpha$ and $\beta$ are separating curves with $\alpha$ an outer curve, (2) both $\alpha$ and $\beta$ are nonouter separating curves, and (3) both $\alpha$ and $\beta$ are nonseparating curves.

\emph{Case 1}: Let $P$ be a pants decomposition for $S$ containing both $\alpha$ and $\beta$. Then $\alpha$ is a vertex of degree 1 in $\mathcal{A}(P)$ which is adjacent only to $\beta$, these properties are preserved by simplicial isomorphism. Therefore $\phi(\alpha)$ is an outer curve adjacent only to $\phi(\beta)$ and hence these curves must form a peripheral pair.

\emph{Case 2}: This is an immediate result of the fact that $\phi$ is an isomorphism and the following technical lemma, which gives a simplicial characterisation of peripheral pairs formed by nonouter separating curves. Recall that the \emph{link} of a vertex $\alpha\in\mathcal{C}^0(S)$ is the complete subgraph of $\mathcal{C}(S)$ induced by all vertices adjacent to $\alpha$ in $\mathcal{C}(S)$. We denote it by $L(\alpha)$. Remark that $L(\alpha)$ is naturally isomorphic to $\mathcal{C}(S \backslash \alpha)$. For any subgraph $\Gamma\subset \mathcal{C}(S)$, we define $\Gamma^*$ as the graph whose vertices are $V(\Gamma)$ and two vertices span an edge if they \textbf{do not} span an edge in $\Gamma$.


\begin{lemma}
	\label{Lemma:PeripheralPairsNonOuter}
Let $\alpha$ and $\beta$ be disjoint nonouter separating curves. Then $(L(\alpha)\cap L(\beta))^*$ has 2 connected components if and only if $\{\alpha,\beta\}$ forms a peripheral pair.
\end{lemma}
\begin{proof}
The necessity of the statement is evident. To prove the sufficiency, remark that $S_{\{\alpha,\beta\}}=S_1\sqcup S_2\sqcup S_3$ . Since 
$(L(\alpha)\cap L(\beta))^*$ has 2 connected components there exists $j$ such that $S_j$ has nonpositive topological complexity. Moreover, given that $\alpha$ and $\beta$ are nonouter separating curves we have that $\alpha\cup\beta=\partial S_j$. A straightforward calculation of the possible topological types for $S_j$ gives us the desired
result.  
\end{proof}

\emph{Case 3}: Up to homeomorphism we can find a separating curve $\gamma$ such that $\{\alpha,\beta,\gamma\}$ bound a pair of pants as in Figure \ref{FIG:NonSepPeriphPair}. Let $P$ be a pants decomposition containing $\{\alpha,\beta,\gamma\}$. By construction and the fact that $\tilde{\phi}$ is a graph isomorphism, we have that $S' \backslash \phi(\gamma)=S_1\sqcup S_2$ and w.l.o.g. we can suppose that $S_1$ has topological complexity equal to 2 and contains $\phi(\alpha)\cup\phi(\beta)$. Since both $\phi(\alpha)$ and $\phi(\beta)$ are nonseparating curves, $S_1$ has positive genus. Therefore $S_{1}$ is homeomorphic to a torus with one boundary component (the boundary curve $\phi(\gamma)$) and one puncture, and the result follows. 
\end{proof}

\section{Proof of Theorem \ref{THM:Main}}\label{Section:Surjectivity}


In this section we use the topological results from the previous section to prove Theorem \ref{THM:Main}. 

Let $\{S_{i}\}$ be a fixed principal exhaustion of $S$. For each $1 \leq i$ we denote by $B_{i}$ the set of boundary curves of $S_{i}$ and $B := \bigcup_{1 \leq i} B_{i}$. 
 
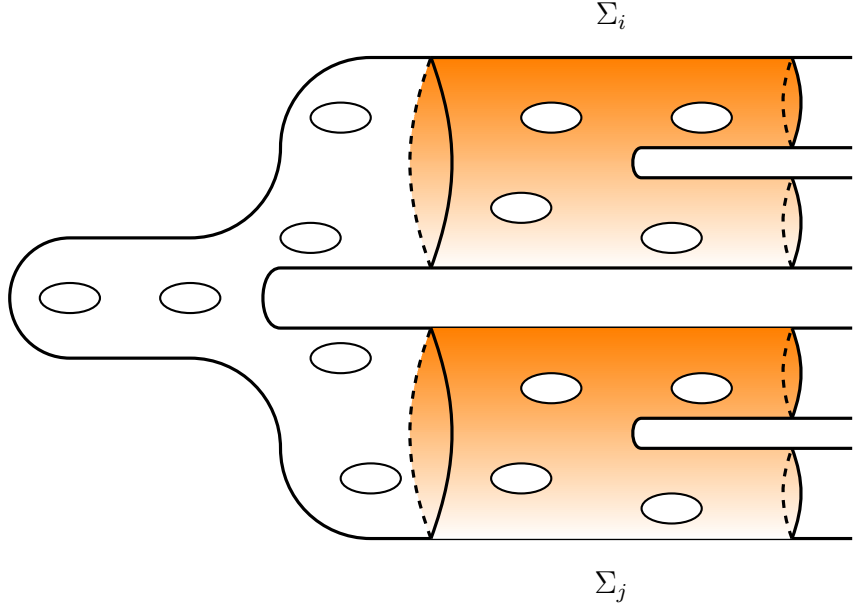
\begin{figure}[!ht]
\begin{center}
\begin{tikzpicture}[scale=.8]

\shadedraw[top color=orange, draw=none] (7,4.5) to[out=110,in=250] (7,8) to[out=0,in=180] (13,8) to[out=290,in=70] (13,6.5) to[out=180,in=0] (10.5,6.5) to[out=180,in=180] (10.5,6) to[out=0,in=180] (13,6) to[out=290,in=70] (13,4.5) to[out=180,in=0] (7,4.5);

\draw[fill=white, thick] (9,7) ellipse (.5 and .25);
\draw[fill=white, thick] (8.5,5.5) ellipse (.5 and .25);
\draw[fill=white, thick] (11.5,7) ellipse (.5 and .25);
\draw[fill=white, thick] (11,5) ellipse (.5 and .25);

\draw[very thick] (13,6.5) to[out=180,in=0] (10.5,6.5) to[out=180,in=180] (10.5,6) to[out=0,in=180] (13,6);
\draw[very thick, dashed] (7,4.5) to[out=110,in=250] (7,8);

\draw[very thick] (7,8) to[out=290,in=70] (7,4.5);

\draw[very thick, dashed] (13,6.5) to[out=110,in=250] (13,8);

\draw[very thick] (13,8) to[out=290,in=70] (13,6.5);

\draw[yshift=-2cm,very thick, dashed] (13,6.5) to[out=110,in=250] (13,8);

\draw[yshift=-2cm,very thick] (13,8) to[out=290,in=70] (13,6.5);

\draw[very thick] (7,8)--(13,8);

\draw[yshift=-3.5cm,very thick] (7,8)--(13,8);

\draw[yshift=-4.5cm,very thick] (7,8)--(13,8);

\draw[yshift=-8cm,very thick] (7,8)--(13,8);


\shadedraw[yshift=-4.5cm,top color=orange, draw=none] (7,4.5) to[out=110,in=250] (7,8) to[out=0,in=180] (13,8) to[out=290,in=70] (13,6.5) to[out=180,in=0] (10.5,6.5) to[out=180,in=180] (10.5,6) to[out=0,in=180] (13,6) to[out=290,in=70] (13,4.5) to[out=180,in=0] (7,4.5);

\draw[yshift=-4.5cm,fill=white, thick] (9,7) ellipse (.5 and .25);
\draw[yshift=-4.5cm,fill=white, thick] (8.5,5.5) ellipse (.5 and .25);
\draw[yshift=-4.5cm,fill=white, thick] (11.5,7) ellipse (.5 and .25);
\draw[yshift=-4.5cm,fill=white, thick] (11,5) ellipse (.5 and .25);

\draw[yshift=-4.5cm,very thick] (13,6.5) to[out=180,in=0] (10.5,6.5) to[out=180,in=180] (10.5,6) to[out=0,in=180] (13,6);

\draw[yshift=-4.5cm,very thick, dashed] (7,4.5) to[out=110,in=250] (7,8);

\draw[yshift=-4.5cm,very thick] (7,8) to[out=290,in=70] (7,4.5);

\draw[yshift=-4.5cm,very thick, dashed] (13,6.5) to[out=110,in=250] (13,8);

\draw[yshift=-4.5cm,very thick] (13,8) to[out=290,in=70] (13,6.5);

\draw[yshift=-6.5cm,very thick, dashed] (13,6.5) to[out=110,in=250] (13,8);

\draw[yshift=-6.5cm,very thick] (13,8) to[out=290,in=70] (13,6.5);


\draw[very thick] (7,0) to[out=180,in=0] (6,0) to[out=180,in=270] (4.5,1.5) to[out=90,in=0] (3,3) to[out=180,in=0] (1,3) to[out=180,in=270] (0,4) to[out=90,in=180] (1,5) to[out=0,in=180] (3,5) to[out=0,in=270] (4.5,6.5) to[out=90,in=180] (6,8) to[out=0,in=180] (7,8);
\draw[very thick] (7,3.5) to[out=180,in=0] (4.5,3.5) to[out=180,in=180] (4.5,4.5) to[out=0,in=180] (7,4.5);

\draw[thick] (1,4) ellipse (.5 and .25);
\draw[thick] (3,4) ellipse (.5 and .25);
\draw[thick] (5.5,7) ellipse (.5 and .25);
\draw[thick] (5,5) ellipse (.5 and .25);
\draw[thick] (5.5,3) ellipse (.5 and .25);
\draw[thick] (6,1) ellipse (.5 and .25);

\draw[very thick] (13,8)--(14,8);
\draw[very thick] (13,6.5)--(14,6.5);
\draw[very thick] (13,6)--(14,6);
\draw[very thick] (13,4.5)--(14,4.5);

\draw[very thick] (13,3.5)--(14,3.5);
\draw[very thick] (13,2)--(14,2);
\draw[very thick] (13,1.5)--(14,1.5);
\draw[very thick] (13,0)--(14,0);

\node at (10,8.7) { $\Sigma_i$};

\node at (10,-.8) { $\Sigma_j$};
 \end{tikzpicture}
\end{center}
\caption{The collection of subsurfaces $\{\Sigma_j\}_{j\geq 1}$.}
\label{FIG:SigmaJotas}
\end{figure}

\begin{theorem}
	\label{THM:Boundaries}
 Let $S$ and $S'$ be infinite-type connected orientable with surfaces with empty boundary and $\phi:\mathcal{C}(S)\to\mathcal{C}(S')$ a simplicial isomorphism. Then $S$ is homeomorphic to $S'$. Moreover, we can construct a homeomorphism $\fun{f}{S}{S'}$ such that $\phi(\beta) = f(\beta)$ for all $\beta \in B$.
\end{theorem}
\begin{proof}
Let $S \backslash B = \bigsqcup_{1 \leq j} \mathrm{int}(\Sigma_{j})$, where the collection $\{\Sigma_{j}\}_{1 \leq j}$ is formed by closed subsurfaces of $S$ of finite type whose topological complexity is at least 4 and such that for any $j \neq k$, $\Sigma_{j} \cap \Sigma_{k}$ is either empty or formed by boundary curves of $S_{i}$, for some $i\in\mathbb{N}$. See Figure \ref{FIG:SigmaJotas}. For each $1 \leq j$, let $P_{j}$ be a pants decomposition of $\Sigma_{j}$ which contains a multicurve analogous to the one used in Proposition \ref{PROP:SameGenus}; in other words we choose a pants decomposition $P_j$ that captures the genus of $\Sigma_j$. Then $P = (\bigcup_{1 \leq j} P_{j}) \cup (\bigcup_{1 \leq i} B_{i})$ is a pants decomposition for $S$.  Let $\mathcal{A}(P)$ be the adjacency graph of $P$ and $\tilde{\phi}:\mathcal{A}(P)\to\mathcal{A}(\phi(P))$ be the corresponding graph isomorphism. 
Curves in $B$ are by definition nonouter and separating, hence every element  $v\in V(B)$ is a cut vertex (\emph{i.e.} $\mathcal{A}(P)\setminus \{v\}$ is disconnected). Then $\mathcal{A}(P)\setminus V(B)= \bigsqcup_{1 \leq j}\Gamma_j$, where each $\Gamma_j$ is a finite subgraph whose vertex set $V(\Gamma_j)$ is precisely the pants decomposition $P_j$ of $\Sigma_j$. By defining $\Sigma_j'$ as the closed subsurface of $S'$ bounded by $\phi(\partial\Sigma_j)$ and recalling that 
 $\phi$ sends nonouter separating curves to nonouter separating curves (see lemma \ref{Lemma:SpecialToSpecial}) we have that
$$
\tilde{\phi}(P)=\bigcup_{1 \leq j} \tilde{\phi}(P_{j}) \cup \bigcup_{1 \leq i} \tilde{\phi}(B_{i})
$$
is such that $\tilde{\phi}(P_j)$ is a pants decomposition for $\Sigma_j'\subset S'$. Since $P_j$ captures the genus of $\Sigma_j$, $\phi$ is an isomorphism and by construction both surfaces have the same number of boundary components, a direct calculation of the topological complexity of $\Sigma_j$ and $\Sigma_j'$ shows that they must be homeomorphic. Moreover, by adjacency w.r.t. $P$, we have that $\partial \Sigma_j'=\{\phi(\alpha):\alpha\subset\partial\Sigma_j\}$. Hence we can find a collection of orientation preserving homeomorphisms $\{\fun{f_{j}}{\Sigma_{j}}{\Sigma_{j}^{\prime}}\}$ such that each $f_{j}$ maps a boundary curve $\alpha\subset\partial\Sigma_{j}$ to $\phi(\alpha)$. These homeomorphisms can be glued together to define a global homeomorphism $f:S\to S'$ which coincides with $\phi$ on $B$.

\end{proof}
 
With this result we have proved the topological rigidity, and we only need to prove that isomorphisms between curve complexes are geometric. 
 
 Hereafter, $\fun{f}{S}{S'}$ denotes the homeomorphism obtained from Theorem \ref{THM:Boundaries}. Remark that  every homeomorphism $h$ of the form $f \circ g$ with $g \in \mathrm{stab}_{pt}(B)$, where
 $$
 \mathrm{stab}_{pt}(B):=\{g\in Homeo(S): \text{$g$ fixes $B$ pointwise}\},
 $$
also coincides with $\phi$ on $B$.
 
For every subsurface $\Sigma$ of $S$ with topological complexity at least $2$, we have that the natural inclusion $\fun{\iota}{\Sigma}{S}$ induces a simplicial map $\fun{\iota_{*}}{\ccomp{\Sigma}}{\ccomp{S}}$ that is an isomorphism on its image. Abusing notation, we denote by $\ccomp{\Sigma}$ the image of $\iota_{*}$ on $\ccomp{S}$. Analogously, we do the same for subsurfaces of $S'$.
\begin{lemma}
\label{Lemma:DentroQuedaDentro}
 For all $1 \leq i$, and for all curves $\alpha \in \ccomp{S_{i}} < \ccomp{S}$, we have that $\phi(\alpha) \in \ccomp{f(S_{i})} < \ccomp{S'}$. In particular, for each $1 \leq i$, the restriction of $\phi$ to $\mathcal{C}(S_i)$
defines an injective simplicial map $\fun{\phi_{i}}{\ccomp{S_{i}}}{\ccomp{f(S_{i})}}$.
\end{lemma}
\begin{proof}
 Let $1 \leq i$ be fixed, $P_{1}$ be a pants decomposition of $S_{i}$, and $P_{2}$ be a pants decomposition of $S \backslash S_{i}$. Then $P = P_{1} \cup B_{i} \cup P_{2}$ is a pants decomposition of $S$.
 
 Recall that $f$ coincides with $\phi$ on $B$, and that the curves in $\partial(f(S_{i}))$ are all separating curves. Then, by the same argument as in Theorem \ref{THM:Boundaries} we have that $\phi(P_{1})$ is a pants decomposition of $f(S_{i})$, and is contained in the interior of $f(S_{i})$. Analogously, the curves in $\phi(P_{2})$ are contained in the interior of $f(S\backslash S_{i})$.
 
 
 Now, let $\alpha$ be a curve contained in $S_{i}$. If $\alpha \in P_{1}$, then $\phi(\alpha) \in \ccomp{f(S_{i})}$ as above. If $\alpha \notin P_{1}$, then there exists $\beta \in P_{1}$ such that $i(\alpha,\beta) \neq 0$. Since we have that:
\begin{itemize}
 \item $\phi(\alpha)$ is disjoint from every element in $\phi(B_{i}) = f(B_{i})$,
 \item $\phi(\beta)$ is contained in $f(S_{i})$, and 
 \item $i(\phi(\alpha),\phi(\beta)) \neq 0$, 
\end{itemize}
we can conclude that $\phi(\alpha)$ is contained in $f(S_{i})$.
\end{proof}
 
 With this lemma and Shackleton's result on combinatorial rigidity (see Theorem 1 in \cite{Shack}), we obtain for each $1 \leq i$ a homeomorphism $\fun{g_{i}}{S_{i}}{f(S_{i})}$  that induces $\phi_{i}$, that is, such that  for all $\alpha\in\mathcal{C}(S_i)$ we have that $\phi(\alpha)=g_i(\alpha)$.
 
 We affirm that each $g_i$ can be extended to a homeomorphism $\overline{g_i}:\overline{S_i}\to \overline{f(S_i)}$ between the closure on $S$ and $S'$ of the respective subsurfaces. To show this we describe first the only possible obstruction to this extension and then why this obstruction never happens. 
 
Each $S_i$ is an open subsurface of $S$ and its punctures can be classified into two categories: those that persist when we take the closure $\overline{S_i}$ of $S_i$ in $S$ (which are precisely those punctures of $S_i$ which are also punctures of $S$) and those that do not (these ``become''  curves contained in $B_i=\partial S_i$ when taking the closure of $S_i$ in $S$). The obstruction could be that $g_i$ exchanges a puncture of $S_i$ that persist in $\overline{S_i}$ with one that does not. We suppose this is the case and we derive a contradiction.

\begin{figure}[!ht]
\begin{center}
\begin{tikzpicture}[scale=2]
\shadedraw[top color=orange, draw=none] (0,0) to[out=90,in=270] (0,.5) to[out=180,in=270] 
(-.5,1.5) to[out=90,in=270] 
(.5,3) to[out=270,in=90] 
(1.5,1.5) to[out=270,in=0] 
(1,.5) to[out=270,in=90] 
(1,0);
\draw[very thick] (0,0) to[out=90,in=270] (0,.5) to[out=180,in=270] 
(-.5,1.5) to[out=90,in=270] 
(.47,3);

\draw[very thick] (.53,3) to[out=270,in=90] 
(1.5,1.5) to[out=270,in=0] 
(1,.5) to[out=270,in=90] 
(1,0);
\draw[fill=white, very thick] (.5,1.5) ellipse (.25 and .5);
\draw (-.5,1.5) to[out=340,in=200] (.25,1.5);
\draw[dashed] (-.5,1.5) to[out=20,in=160] (.25,1.5);
\draw[xshift=1.25cm] (-.5,1.5) to[out=340,in=200] (.25,1.5);
\draw[xshift=1.25cm, dashed] (-.5,1.5) to[out=20,in=160] (.25,1.5);
\draw (0,.5) to[out=340,in=200] (1,.5);
\draw[dashed] (0,.5) to[out=20,in=160] (1,.5);

\node at (-0.75,1.5) {$\beta$};
\node at (1.75,1.5) {$\gamma$};
\node at (-0.25,.25) {$\alpha$};
\draw[fill=white, thick] (.5,3) ellipse (.05 and .025);
\end{tikzpicture}
\end{center}
\caption{$\{\beta,\gamma\}$ is a peripheral pair.}
\label{FIG:PeriphPairBis}
\end{figure}
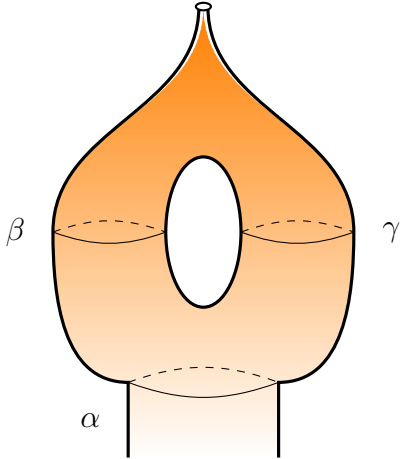

Let $\alpha,\beta,\gamma$  bound a pair of pants in $S$ such that $\alpha\subset\partial S_i$ and $\beta,\gamma\in\mathcal{C}(S_i)$. Note this implies that $\{\beta,\gamma\}$ is a peripheral pair in $S_{i}$. If $g_i$ exchanges the puncture of $S_i$ defined by $S_i\setminus \alpha$ with a puncture of $S$, then $\{\phi(\beta),\phi(\gamma)\}$ would be a peripheral pair in $S$. By the Proposition \ref{FIG:NonSepPeriphPair}, $\{\beta,\gamma\}$ is also a peripheral pair in $S$. 
This situation is depicted in Figure \ref{FIG:PeriphPairBis}. It is clear that $S\setminus\alpha$ has one connected component whose topological complexity is \textit{strictly less} than 3. This is a contradiction, for both connected components of $S\setminus\alpha$ have topological complexity \textit{at least} 3. 

Thus for each $1 \leq i$, we have a homeomorphism $\fun{\overline{g}_{i}}{\overline{S_{i}}}{\overline{f(S_{i})}}$ that induces $\phi_{i}$. Using the following lemma we can assert that for each $1\leq i$, $\overline{g_i}$ coincides with $f$ on $B_i$.
 
\begin{figure}[!ht]
\begin{center}
\begin{tikzpicture}[scale=1]

\shadedraw[top color= orange, very thick]
(0,0) to[out=45,in=270] (.5,1) to[out=90,in=180] (1.5,4)
to[out=340,in=200] (2.5,4)
to[out=0,in=90] (3.5,1) to[out=270,in=135] (4,0) to[out=135,in=45] (3,0)
to[out=90,in=0] (2,1)
to[out=180,in=90] (1,0) to[out=135,in=45] (0,0)
;
\draw[very thick] (0,0) to[out=315,in=225] (1,0);
\draw[xshift=3cm,very thick] (0,0) to[out=315,in=225] (1,0);
\draw[very thick] (1.5,4)
to[out=20,in=160] (2.5,4);

\draw[fill=white, thick] (2,2) to [out=180,in=270] (1.5,2.5) to [out=90,in=180] (2,3)
to [out=0,in=90] (2.5,2.5) to [out=270,in=0] (2,2); 

\draw[thick] (1.5,2.5) to [out=270,in=45] (.5,1);
\draw[thick, dashed] (1.5,2.5) to [out=225,in=70] (.5,1);

\draw[thick] (2,2) to [out=315,in=45] (2,1);
\draw[thick, dashed] (2,2) to [out=225,in=135] (2,1);

\draw[thick] (2.5,2.5) to [out=300,in=160] (3.5,1);
\draw[thick, dashed] (2.5,2.5) to [out=340,in=130] (3.5,1);

\node at (0,1.5) {$\phi(\gamma)$};

\node at (2,.5) {$\phi(\alpha)$};

\node at (4,1.5) {$\phi(\beta)$};

 \end{tikzpicture}
\end{center}
\caption{ the curves $\phi(\alpha)$, $\phi(\beta)$ and $\phi(\gamma)$ do \textit{not} bound a pair of pants.}
\label{FIG:NotBoundingPairofPants}
\end{figure}
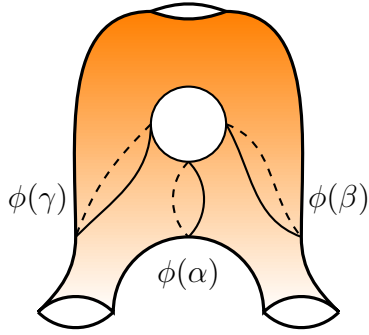

 \begin{lemma}
 	\label{Lemma:BoundPants}
 Let $\alpha$, $\beta$ and $\gamma$ be curves on $S$ such that $\alpha$ is a separating curve, and $\alpha$, $\beta$ and $\gamma$ bound a pair of pants on $S$. Then $\phi(\alpha)$, $\phi(\beta)$ and $\phi(\gamma)$ also bound a pair of pants on $S'$.
\end{lemma}
\begin{proof}
 Let $P$ be a pants decomposition of $S$ with $\alpha, \beta, \gamma \in P$. Then, with respect to $P$, $\alpha$ is adjacent to $\beta$, $\beta$ is adjacent to $\gamma$, and $\alpha$ is adjacent to $\gamma$. By 
 Proposition \ref{PROP:GraphIsomorph}, we know that adjacency is preserved under $\phi$.
 
 The only possibility for this to happen and having that $\phi(\alpha)$, $\phi(\beta)$ and $\phi(\gamma)$ do \textit{not} bound a pair of pants on $S'$, is (up to homeomorphism) illustrated in Figure \ref{FIG:NotBoundingPairofPants}. However, if this were to happen we could find a curve $\delta$ on $S'$ that would intersect $\phi(\alpha)$ exactly once, which is impossible since $\phi(\alpha)$ is a separating curve due to Proposition \ref{PROP:GraphIsomorph}. Therefore, $\phi(\alpha)$, $\phi(\beta)$ and $\phi(\gamma)$ bound a pair of pants on $S'$.
\end{proof}


For each $1\leq i$ we can pick an adequate element $l_i$ of the form $f \circ \eta$ with $\eta \in \mathrm{stab}_{pt}(B)$ and define:
$$
h_i(x)=\begin{cases}
\overline{g_i}(x) &\text{ if $x\in S_i$}\\
l_i(x) &\text{ if $x\in S\setminus S_i$}
\end{cases}
$$
We obtain this way a family of homeomorphisms $h_i:S\to S'$ which by construction satisfy that $h_i(\alpha)=\phi(\alpha)=h_j(\alpha)$ for all $i<j$ and $\alpha\in\mathcal{C}(S_i)\subset\mathcal{C}(S_j)\subset\mathcal{C}(S)$. As a consequence of the Alexander method (see \cite{FarbMar}, chapter 2.3) we have for each $i<j$ that $h_{i\vert S_i}=h_{j\vert S_j}\circ M_{i}$, where $M_i\in Homeo(S_i)$ is a multitwist whose support is contained in a neighbourhood in $S_i$ of $\partial S_i$. In other words, for each $1\leq i$ there exists a subsurface $\widetilde{S_i}\subset S_i\subset S$ isotopic within $S_i$ to $S_i$ such that the support of the multitwist $M_i$ is contained in $S\setminus\widetilde{S_i}=\sqcup_{k=1}^s A_k$, where each $A_k$ is an annulus. In particular $M_{\vert \widetilde{S_i}}=Id_{\widetilde{S_i}}$ and hence for each $i<j$ we have that $h_{i\vert\widetilde{S_i}}=h_{j\vert\widetilde{S_i}}$.  This way we can define the following map:
\begin{center}
 \begin{tabular}{rccc}
  $h:$ & $S$ & $\longrightarrow$ & $S'$\\
   & $s \in \widetilde{S}_{i}$ & $\mapsto$ & $h_{i}(s)$
 \end{tabular}
\end{center}
 Since for all $1 \leq i < j$ we have that $h_{i}|_{\widetilde{S}_{i}} = h_{j}|_{\widetilde{S}_{i}}$, this map is well-defined. Moreover it is a homeomorphism and by construction it coincides with $\phi$ on the whole curve graph $\ccomp{S}$, as desired.
\qed

\bibliographystyle{plain}
\bibliography{biblio}

\end{document}